\documentclass[11pt]{amsart}
\usepackage{amsmath,amssymb,amscd}
\usepackage{amsmath,amssymb}

\usepackage{color}
\usepackage{amsxtra, amsmath}
\usepackage{amssymb, amscd}
\usepackage{graphicx}
\usepackage{mathrsfs}

\newcommand\NN{\mathbb N}

\theoremstyle{plain}
\newtheorem{thm}{Theorem}[section]

\theoremstyle{definition}

\theoremstyle{remark}
\newtheorem{remark}[thm]{Remark}

\title[A new commutativity property of exceptional orthogonal polynomials ]{
A new commutativity property of exceptional orthogonal polynomials}

\author{M.M. Castro$^*$$^1$ and  F. A.  Gr\"unbaum$^2$}



\begin{document}

\begin{abstract} We exhibit three examples showing that the ``time-and-band limiting" commutative property found and exploited
	by  D. Slepian, H. Landau and H. Pollak at Bell Labs in the 1960's,  and independently by M. Mehta and later by C. Tracy and H. Widom in Random matrix theory, holds for exceptional orthogonal polynomials. The property in question is the existence of local operators with simple spectrum that commute with naturally appearing global ones.
	We illustrate numerically the advantage of having such a local operator. 
\end{abstract}

\maketitle

\bigskip

\begin{center}
$^1$$ ^*$Departamento de Matem\'atica Aplicada II and IMUS, Escuela Polit\'ecnica Superior,\\ 
 Universidad de Sevilla, c. Virgen de \'Africa 7, 41011, Seville, Spain.

\bigskip

$^2$Department of Mathematics, University of California, Berkeley
	CA 94705.

\bigskip

 $^*$Corresponding author E-mail: mirta@us.es;

\medskip

$^2$grunbaum@math.berkeley.edu;

\end{center}

\bigskip

\section{A brief historical introduction}\label{intro}

The interest of one of us in what is now called the ``bispectral problem"
posed and solved in \cite{DG} arose from an effort to understand and 
extend a mathematical miracle, uncovered by D. Slepian, H. Landau and 
H. Pollak at Bell Labs back in the 1960's,  \cite{SLP2,SLP3,SLP4,S2,SLP5,S1,SLP1}.  It is of obvious importance in signal processing and it was motivated by
work of C. Shannon, see \cite{Sh}.
This is mentioned in the introduction to \cite{DG} and is recalled in the next few lines.

\bigskip

For an unknown
signal $f(t)$ supported in $[-T,T]$ one observes its Fourier transform
$Ff(k)$ for frequencies $k$ in the band $[-\mathcal{W},\mathcal{W}]$. The numerically stable
reconstruction of $f$ from this data leads to the study of an integral operator in
$L^2\left(-T,T\right)$ with kernel given by $K(t,s)=\sin \mathcal{W}(t-s)/(t-s)$. One needs to
compute numerically many of its eigenfunctions, and this gives a very
ill-conditioned problem because the eigenvalues are (except for
about $4T\mathcal{W}$ of them) very close together. The miracle in question is that
one can exhibit a selfadjoint second order differential operator that
commutes with the integral operator, has a simple and well spread out spectrum resulting in a 
common set of
eigenfunctions. One has replaced a very ill-posed problem by a very
well posed one.

\bigskip

The topic of exceptional orthogonal polynomials has developed into an
active area of research in the last few years. The literature is large
and we just mention a few papers, \cite{Durec,Du,GGU2,GUHerm,GU2,GUKKM,KM,Qu,STZ}.
The {\it potential} applications of these polynomials are quite varied: they
provide a vast extension of the classical orthogonal polynomials of
Jacobi, Laguerre and Hermite which feature in countless areas of mathematics,
both pure and applied. They each give a basis of polynomials that are joint
eigenfunctions of a fixed differential operator $L$ of order two

$$L p_{n}(x) = \lambda_n  p_{n} (x)$$

The main difference with the classical ones is that the index $n$, that 
indicates the degree of $p_{n}$,  needs not run over the entire set
${0,1,2,3,...}$. A way to increase the chances that this larger class of polynomials will be widely used is to study
 the extend to which they share certain properties with their classical counterparts.
Such an effort will necessarily develop in
an {\it exploratory} fashion. This is the spirit of this paper.

\bigskip

The hope at the time of \cite{DG} was that situations exhibiting the highly unusual ``bispectral property", to
be defined below, would give rise to extensions of the commutativity miracle mentioned above beyond the Fourier case exploited by the Bell Labs group.
One says that one has a ``bispectral situation" when a differential operator
$L$ has eigenfunctions $f(x,k)$ that satisfy a differential equation in
the spectral parameter $k$, i.e. we have

\begin{equation}\label{bispL}
	L f(x,k) = \Lambda(k) f(x,k),
	\end{equation}

as well as
\begin{equation}\label{bispB}
 B f(x,k) = \Theta(x) f(x,k) ,
\end{equation}
where $B$ is a differential or difference operator acting on $k$. The spectral parameter $k$ runs over a discrete or
continuous set.

\bigskip

Several tools were used in \cite{DG} to classify all the
situations when $L$ has order two and $B$ is a differential operator of arbitrary order. They include, among other things, the so called Darboux
processes, the so called ad-conditions, as well as a careful study of the monodromy properties of $L$.
The main surprise was the observation that the rational solutions of the
Korteweg-deVries equation play a important role in half of the cases found in \cite{DG}.
The role of its master symmetries was observed in a later paper \cite{ZM}.

The case when $k$ is a discrete variable and B is a second order difference operator, as in the case of the classical orthogonal polynomials, was considered in a series of
papers \cite{GH1,GH2,GHH4,GHH3}, where the Darboux process was applied to
either the  semi-infinite or doubly infinite banded matrix $B$  to obtain the Krall polynomials or functions that extend them. Here the
role of the KdV equation is taken up by nonlinear evolutions such as the Toda flows. The appearance of these integrable isospectral systems when $B$ is a differential or a difference operator came about in these papers by using the full power of the Darboux process (see also \cite{HP,PTesis}). Two earlier papers by M. Reach \cite{R1,R}, based on his UC Berkeley thesis $1987$, deal with applying the Darboux process
to the second order operator $L$, as in \cite{DG}, but allowing $B$ to be a difference
recursion of arbitrary order. Just as in \cite{DG} this leads to increasing orders in the
recursion relation given by $B$. See also \cite{GUKKM,KM} where one has a similar situation.

\bigskip

This is not the place to review in detail the developments just mentioned, and we just recall the 
very important contributions by G. Wilson, \cite{W1,W2}, A. Kasman
and M. Rothstein, \cite{KR}, as well as those of B. Bakalov, E. Horozov and M. Yakimov,
\cite{BHY1,BHY2,BHY3}. More references can be found in \cite{G8,G6,HK}.

It is clear that a certain set of common tools were used both in the study of the bispectral
problem as well as in the study of exceptional orthogonal polynomials.
We just mention a few of these.
Apparently the first explicit mention of the bispectral
property in the second area is made in \cite{STZ} (see section $4$ as
well as the acknowledgments in this reference). It appears that the first paper
to exploit the Darboux process in connection with exceptional orthogonal 
polynomials is \cite{Qu}.
The role of trivial monodromy has appeared in later papers such as \cite{GGU2,GUHerm}. One can see points of contact between the considerations in
\cite{GUKKM,KM} and those in the bispectral problem as mentioned by these authors (see also \cite{GUHerm}).
There are many differences between these two topics: all considerations
in the bispectral problem are of a local nature, whereas in the case of exceptional
orthogonal polynomials the issue of the completeness of these polynomials in some
appropriate Hilbert space is of importance.

\bigskip

Making heavy use of the bispectral property, the paper \cite{G3} establishes
that the commutativity phenomenon alluded to above holds for the classical orthogonal
polynomials. The paper \cite{G4} shows that the same property holds in
connection with the ``even family" in \cite{DG}, the one connected with
the master symmetries of the KdV equation. A general strategy to connect
bispectrality and the commutativity property
is given in \cite{GY}. A much more recent push
in this direction is in \cite{CGYZ1,CGYZ2,CGYZ3,CGYZ4}.

\bigskip

The trivial cases of bispectrality when both physical and frequency space are the real line are given by  Bessel and Airy operators (the Bessel case includes Fourier analysis). They have featured in important problems of mathematical physics, such as potential theory, electromagnetism and optics for a very long time. Both cases lead to integral operators admitting a commuting differential operator. For the Bessel case this was proved by D. Slepian, see \cite{SLP5}. The Airy case was observed by C. Tracy and H. Widom in the context of Random matrix theory, see \cite{TW1}. In \cite{CGYZ3} one considers deformations of the Airy integral operator which preserve the commutativity property. For the Bessel case see \cite{G4}. For very recent numerical work on the eigenfunctions of the Airy integral operator which exploits the existence of the commuting one see \cite{ShS}. This work has applications both in Random matrix theory as well as in optics.


\bigskip

In view of these more recent papers, and once
one notices that exceptional orthogonal polynomials give a bispectral situation and are connected to the classical ones by applications of the Darboux process the
results in the present paper are not totally unexpected. However, the nature of the explicit results given here
makes it worthwhile presenting them separately. One should mention that there are matrix valued versions of the commutativity property   too,
see \cite{CG5, CG6, GPZ1,GPZ2,GPZ3}.

\bigskip

A last historical remark: many ideas appear again and again in the development of any area of mathematics and the issue of giving appropriate credit 
takes complicated turns. It is nice to be able to single out a paper by V.
Bargmann, \cite{B}, where several of the concerns and issues in modern spectral theory originate. The very nice paper by C. Quesne
\cite{Qu} starts by referring to that pioneering paper. In fact, the introduction to her paper offers a rather instructive view of a rich area where many tools appear over and over again. A very complete discussion of the work of V.Bargmann, 
M.G. Krein, V.A. Marchenko, I.M. Gelfand and B.M. Levitan (among others) is given in the book \cite{ChS}.

\bigskip

In conclusion, we mention that the phenomenon mentioned
above (in its original Fourier version) has found a rather unexpected use in a series of papers by A. Connes    and collaborators in connection with the Riemann zeta function.
For the most recent push in this direction, see \cite{CM}.
For a commentary on this paper see \cite{G9}. For a recent use of this commutativity property (again in its original
 Fourier form) in connection with the Bethe ansatz and entanglement, see \cite{BP,BCV,CNV} and its
references. It is interesting that many of the kernels that appear in connection with the bispectral problem play an important role in Random Matrix Theory, see \cite {TW1,TW2}. For deformations of these kernels that preserve both phenomena see \cite{CGYZ1,CGYZ2,CGYZ3,CGYZ4}. These observations argue for 
approaching the relatively new area of exceptional orthogonal polynomials
by casting a wide net.

\section{The contents of the paper}\label{structure}

The commutativity phenomenon alluded to above  was extended to the case of classical orthogonal polynomials
defined on a finite set, see \cite{P1,P2}. It was later observed by R. Perline that the explicit
form of the commuting operator could be given a simple and unified form, see \cite{Per}.

It was subsequently seen, see \cite{GPZ1, GVZ} that this simple form of Perline applies unchanged to
other situations. In \cite{GVZ} one observes that (as mentioned by R. Perline) the simple
form of the commuting operator may result in one that does not have simple spectrum. This happens
for the so called Bannai-Ito polynomials. In \cite{GVZ} one shows that a more complicated form
of the commuting operator, still built with an appropriate extension of  Perline's construction, yields one with simple spectrum.

Since the examples of exceptional orthogonal polynomials that we consider here, namely Jacobi, Laguerre and 
Hermite, involve recursion relations of orders five and seven respectively, the simple form
put forward by R. Perline has to be modified. This point will be illustrated in some of the examples below.


\bigskip

We are now in a position to describe the contents of the present paper.
In section \ref{preliminars} we recall the definitions of the operator of time-band-time
limiting as well as the operator of band-time-band limiting.  
In section \ref{Hermite} we consider the case of the exceptional Hermite polynomials
which do not depend on free parameters. In section \ref{Jacobi} we consider the simplest instance of the exceptional Jacobi polynomials. These ones depend on the usual
parameters $\alpha$ and $\beta$ and some of our results are given
in terms of these parameters. Some results, which have been checked for multiple
values of the parameters,
are illustrated by a specific (but arbitrary) choice of them. In section \ref{Laguerre}
we consider an instance of exceptional Laguerre polynomials. Finally, in section \ref{numeric}, in the spirit of \cite[Section 6]{CG5}, we conclude by exploiting the numerical pay-off of the results of the previous sections.

\section{The operators of time and band limiting}\label{preliminars}


We start with a very general setup of the time-band limiting problem for orthogonal polynomials (see for instance \cite{G3}) which will be applied later on to the different examples of sequences of exceptional orthogonal polynomials discussed in this paper.

\medskip

Let $w=w(x)$ be a  weight function  in the open interval $(a,b)$, for which all the moments $\displaystyle \int_{a}^bx^nw(x)dx$, $n\geq 0$, are finite. Let $\left(p_n(x)\right)_{n\geq0}$ be a sequence of real valued orthonormal polynomials with respect to the weight $w(x)$. Since we will be dealing with exceptional orthogonal polynomials, we do not assume that $\deg  p_n=n$. Consider the following two Hilbert spaces: The space $L^2((a,b), w(x)dx)$, denoted here by $L^2(w)$, of all measurable functions $f(x)$, $x\in (a,b)$, satisfying $\int_a^b f^2(x)w(x)dx < \infty $ and the space $\ell^2(\NN_0)$
of all real valued sequences $(c_n)_{n\in \NN_0}$ such that $\sum_{n=0}^\infty c_n^2 < \infty$.


The map $\mathcal F:\ell^2(\NN_0) \longrightarrow L^2(w)$ given by
$$(c_n)_{n=0}^\infty \longmapsto  \sum_{n=0}^\infty c_n p_n(x)$$
is an isometry. If the polynomials are dense in $L^2(w)$, this map is unitary with the inverse $\mathcal F^{-1}: L^2(w)\longrightarrow \ell^2(\NN_0) $ given by
$$ f(x) \longmapsto c_n=\int_a^b f(x)p_n(x)w(x) dx.$$



\bigskip

We denote our map by $\mathcal F$ to remind the reader of the usual Fourier transform. Here $\NN_0$ takes up the role of ``frequency space" and the interval $(a,b)$ the role of ``physical space".

\bigskip

The {\em band limiting  operator}, at level $N$ acts on $\ell^2(\NN_0)$ by simply setting equal to zero all the components with index larger than $N$. We denote it by $\chi_N$. The {\em time limiting operator }, at level $T$,  acts on $L^2(w)$ by multiplication by the characteristic function of the interval $(a,T]$, $T\le b$. This operator will be denoted by $\chi_T$.

\medskip

Consider the problem of determining a function $f$  from the following data: $f$ has support on the finite set $\{0,\dots , N\}$ and its Fourier transform $\mathcal Ff$ is known on the  set $(a,T]$. This can be formalized as follows
$$\chi_T \mathcal Ff=g=\text{ known},\qquad  \chi_N f=f.$$

We can combine the two equations into

$$E f= \chi_T \mathcal F \chi_N f=g.$$

To analyze this problem we need to compute the singular vectors (and values) of the operator $E:\ell^2(\NN_0)\longrightarrow L^2(w) $. These are given by the eigenvectors of the operators
$$E^*E= \chi_N \mathcal F^{-1} \chi_T \mathcal F \chi_N\qquad \text{ and } \qquad S_2=E E^*= \chi_T \mathcal F \chi_N \mathcal F^{-1} \chi_T.$$

The operator $E^*E$, acting in $\ell^2(\NN_0)$ is just a finite dimensional matrix $M$, and each entry is given by
$$(M)_{m,n}=(E^*E)_{m,n}= \int_a^T p_m(x) p_n(x) w(x) dx, \qquad 0\leq m,n \leq N.$$
The second operator $S_2= E E^*$ acts in $L^2((a,T),w(x)dx)$ by means of the integral kernel
$$k(x,y)=\sum_{j=0}^N p_j(x)p_j(y).$$

Consider now the problem of finding the eigenfunctions of $E^* E$ and/or $E E^*$. For arbitrary $N$ and $T$ there is no hope of doing this analytically, and one has to resort to numerical methods. This is a remarkably ill-conditioned problem since most of the eigenvalues are crowded together. Of all the strategies one can dream of for handling this problem, none sounds so appealing as that of finding a differential operator with simple-and spread out-spectrum which would have the same eigenfunctions as the original operators. This is exactly what Mehta as well as Slepian, Landau and Pollak did  when dealing with the real line and the actual Fourier transform. They discovered (the analog of) the following properties:

\begin{itemize}
	
	\item For each $N$, $T$ there exists a symmetric matrix $\mathscr{L}$ with a small number of diagonals, with simple spectrum, commuting with $M$.
	
	\item For each $N$, $T$ there exists a selfadjoint differential operator $D$, with simple spectrum, commuting with the integral operator $S_2=EE^*$.
	
\end{itemize}

\bigskip

In this paper we will see instances of exceptional orthogonal polynomials where this phenomenon holds. Once more we will see that the ``bispectral property", first considered in \cite{DG}, guarantees the commutativity of  these two operators, a global and a local one.

\medskip

For an up-to-date treatment of the important issue of computing numerically the eigenfunctions of $D$, see \cite{ORokX}. For the case of the DFT, see \cite{G0}.



\section{The exceptional Hermite  polynomials}\label{Hermite}

We consider the family of exceptional Hermite polynomials defined by $\hat{H}_0=1$,
$$\hat{H}_n=H_n+4nH_{n-2}+4n(n-3)H_{n-4}, \quad n\geq 3,$$
where $H_n$ are the classical Hermite polynomials given by the Rodrigues formula

$$H_n=(-1)^ne^{x^2}D_x^ne^{-x^2}, n=0,1,2\ldots .$$

\bigskip

These polynomials can also be defined by means of a determinant (see for instance \cite{GUHerm}):
$$\hat{H}_n= \frac{1}{8(n-1)(n-2)}\textrm{det}\left(\begin{array}{ccc}H_n&H'_n&H''_n \\H_1&H'_1&H''_1\\H_2&H'_2&H''_2 \end{array} \right),\quad n\neq 1,2.$$

The exceptional Hermite polynomials satisfy the orthogonality relation:
$$ \int_{-\infty}^{\infty}\hat{H}_m(x)\hat{H}_n(x)\dfrac{e^{-x^2}}{(1+2x^2)^2}{dx}
=\dfrac{\sqrt{\pi}2^nn!}{(n-1)(n-2)}\delta_{n,m},\quad n\neq 1,2.$$


\subsection{Bispectral operators}


Let us now consider the orthonormal sequence of exceptional Hermite orthogonal polynomials given by $\widetilde{H}_n=\dfrac{\sqrt{(n - 1) (n - 2)}}{\sqrt[4]{\pi}\sqrt{ 2^nn!}}\hat{H}_n$.

The exceptional polynomials $\widetilde{H}_n$ satisfy the  differential equation
\begin{equation*}\label{eqdifHerm}
\widetilde{H}_n''-\left(2x+\dfrac{8x}{1+2x^2}\right)\widetilde{H}_n' +2n\widetilde{H}_n=0
\end{equation*}
and the  recurrence relation, written explicitly in \cite{Durec} for a different normalization:
\begin{equation}\label{recurr}
\alpha_{n-3}\widetilde{H}_{n-3}+\beta_{n-1}\widetilde{H}_{n-1}+\beta_n\widetilde{H}_{n+1}+\alpha_n\widetilde{H}_{n+3}=\Theta(x)\widetilde{H}_n,\quad n\geq 0,
\end{equation}

where the coefficients $\alpha_n$ and $\beta_n$ and the function $\Theta(x)$ are given by:
\begin{equation}\label{coefs}
\alpha_n=\dfrac{\sqrt{2}}{3}\sqrt{(n+3)(n-1)(n-2)},\quad \beta_n=\sqrt{2(n+1)n(n-2)},\quad \Theta(x)= \frac{4x^3}{3} + 2 x .
\end{equation}

Here, we understand that $\widetilde{H}_1(x)=\widetilde{H}_2(x)=0$ as well as $\widetilde{H}_{-2}(x)=\widetilde{H}_{-1}(x)=0$.

\medskip

Notice that $\Theta^{\prime}(x)=2(2x^2+1)=\dfrac{1}{2}\,\det\begin{pmatrix} H_1&H_1'\\H_2&H_2'\end{pmatrix}$. This allows, in the spirit of \cite{DG}, to relate $\Theta(x)$ to the appropriate Sato's $\tau$ function.

\medskip

If we write $$\Psi(x)=(\widetilde{H}_0(x),\widetilde{H}_3(x),\widetilde{H}_4(x),\ldots)^t$$
and B for the semi-infinite heptadiagonal symmetric matrix:

\begin{equation}\label{matrizB}
B=%
\left(
\begin{array}{cccccccccccc}
0&\alpha_0&0&\ldots&0&\ldots& & &&&&\\
\alpha_0 &0&\beta_3&0&\alpha_3&0&\ldots&&&&&\\
0&\beta_3&0&\beta_4&0&\alpha_4&0&\ldots&&&&\\
0&0&\beta_4&0&\beta_5&0&\alpha_5&0&\ldots&&&\\ 0&\alpha_3&0&\beta_5&0& \beta_6&0&\alpha_6&0\ldots&&& \\ 0&0&\alpha_4&0&\beta_6&0&\beta_7&0&\alpha_7&0&\ldots&\\
\vdots&&\ddots&&\ddots& \ddots&&&&&\ddots&
\end{array}
\right),
\end{equation}  

one has $B\Psi(x)=\Theta(x)\Psi(x)$ (see (\ref{bispB})), where the expression of $\Theta(x)$ is given in (\ref{coefs}).

From the other side one has $L\Psi(x)=\Lambda(n)\Psi(x)$, see (\ref{bispL}), where 
\begin{equation}\label{opLHermite}
L=\dfrac{d^2}{dx^2}-\left(2x+\dfrac{8x}{1+2x^2}\right)\dfrac{d}{dx}
\end{equation}

and  $\Lambda(n)$ is the diagonal matrix whose entries are $\lambda_n=-2n$, $n=0,3,4,\ldots $.


\subsection{Time and Band limiting}

Here, we define the $N\times N$ matrix $M$ of truncated inner products depending on a real parameter $T$, whose entries are given by:

$$M_{m,n}=\int_{-\infty}^{T}\widetilde{H}_m(x)\widetilde{H}_n(x)\dfrac{e^{-x^2}}{(1+2x^2)^2}dx,\quad n,m=0,3,4,\ldots, N+1.$$

For  fixed values of the parameters $N$ and $T$ one looks for a ``narrow banded" commuting matrix $\mathscr{L}$.

One finds that there exists an {\it heptadiagonal} matrix $\mathscr{L}$ commuting with $M$. If we consider for instance the normalization $\mathscr{L}_{N,N}=0$ and  $\mathscr{L}_{N,N-3}=1$, this matrix is unique.

\bigskip

We display here the  symmetric time-band limiting matrix $M$ of size 7 that 
appears in the Hermite case. Here and in what follows, $I$ will denote the identity  matrix of appropriate size.

\bigskip
$
M-\left(\dfrac{1+Erf(T)}{2}\right)I=
$


$$
\left(\begin{array}{cccccc}-T& {1}\over{\sqrt{3}} & {{T}}\over{{\sqrt{2}}} &
\frac{2T^2-1}{	\sqrt{10}}&\\ 
{1}\over{\sqrt{3}}& \frac{T(2T^2+3)}{3}&\scriptstyle{\left(4\,T^4+12
	\,T^2+3\right)}\over{2\sqrt{6}}&{{\sqrt{2}\,T^3\,\left(2\,T^2
		+5\right)}\over{\sqrt{15}}}
\\ {{T}}\over{{\sqrt{2}}}&\scriptstyle{\left(4\,T^4+12
	\,T^2+3\right)}\over{2\sqrt{6}} &\scriptstyle{\frac{1}{4} T \left(4 T^4+12 T^2+3\right)}&{{\left(8\,T^6+20\,T^4+10\,T^2+5\right)}\over{
		4\,\sqrt{5}}}\\\frac{2T^2-1}{	\sqrt{10}}&{{\sqrt{2}\,T^3\,\left(2\,T^2
		+5\right)}\over{\sqrt{15}}}&{{\left(8\,T^6+20\,T^4+10\,T^2+5\right)}\over{
		4\,\sqrt{5}}} &\scriptstyle{\frac{T}{5} (5 + 5 T^2 + 8 T^4 + 4 T^6)}
\\{T\,\left(2\,T^2-3\right)}\over{3\,\sqrt{2}\,
}&{ \scriptstyle{8\,T^6+12\,T^4-18\,T^2-3}}\over{
	6\sqrt{6}}&{{T\,\left(8\,T^6+12\,T^4-6\,T^2+3\right)
}}\over{12} & {{16\,T^8+16
		\,T^6+60\,T^2+15}\over{12\,\sqrt{5}}} &
\\{{\left(4\,T^4-12\,T^2+3
		\right)}}\over{{2\,\sqrt{42}}}&\scriptstyle{T^3\,\left(4\,T^4-21\right)}\over{3\,
	\sqrt{14}}&{{\left(16\,T^8-56\,T^4-7\right)}\over{8\,
		\sqrt{21}}} & {{T^3\,
		\left(8\,T^6-4\,T^4-14\,T^2+35\right)}\over{2\,\sqrt{15}
		\sqrt{7}}}  \\{T\,\left(\scriptsize{4\,T^4-20\,T^2+15}
	\right)}\over{{4\,\sqrt{30}}}&\scriptstyle{\left(
	16\,T^8-32\,T^6-120\,T^4+72\,T^2+9\right)}\over{24
	\sqrt{10}}&{{T\,\left(16\,T^8-32\,T^6-88\,T^4+40\,T^2-15\right)}}\over{16\,\sqrt{15}} & {{32\,T^{10}-80\,T^8-80\,T^6+200\,T^4-150\,T^2-25
}}\over{80\,\sqrt{3}}
\end{array} \right.
$$

$$
\begin{array}{cc}
{T\,\left(2\,T^2-3\right)}\over{3\,\sqrt{2}\,
}&{{4\,T^4-12\,T^2+3
}}\over{{2\,\sqrt{42}}}\\\scriptstyle{8\,T^6+12\,T^4-18\,T^2-3}\over{
	6\sqrt{6}}&\scriptstyle{T^3\,\left(4\,T^4-21\right)}\over{3\,
	\sqrt{14}}\\{{T\,\left(8\,T^6+12\,T^4-6\,T^2+3\right)
}}\over{12}&{{16\,T^8-56\,T^4-7}\over{8\,
		\sqrt{21}}}\\
{{16\,T^8+16
		\,T^6+60\,T^2+15}\over{12\,\sqrt{5}}}&{{T^3\,
		\left(8\,T^6-4\,T^4-14\,T^2+35\right)}\over{2\,\sqrt{15}
		\sqrt{7}}}  \\ \scriptstyle{\frac{T}{36}  \left(16 T^8+120 T^2+27\right)}& {{32\,T^{10}-48\,T^8+336\,T^4+126\,T^2+63}}
\over{24\sqrt{21}} \\{{
		32\,T^{10}-48\,T^8+336\,T^4+126\,T^2+63
	}\over{24\,\sqrt{21}}}& \frac{T \scriptstyle{\left(16 T^{10}-48 T^8+40 T^6+252 T^4+77 T^2+84\right)}}{84} \\
{{\scriptstyle {T\,\left(32\,T^{10}-112\,T^8+48\,T^6+408\,T^4-
			150\,T^2+45\right)}}\over{48\sqrt{15
}}}&{{\scriptstyle 64\,T^{12}-320\,T^{10}+464\,T^8+1120\,T^6-420\,T^4+1260\,T^2+
		315}}\over{96\,\sqrt{35}}
\end{array}
$$


$$
\left.\begin{array}{c}
{T\,\left(\scriptsize{4\,T^4-20\,T^2+15}
	\right)}\over{{4\,\sqrt{30}}}\\\scriptstyle{\left(
	16\,T^8-32\,T^6-120\,T^4+72\,T^2+9\right)}\over{24
	\sqrt{10}}\\{{T\,\left(16\,T^8-32\,T^6-88\,T^4+40\,T^2-15\right)}}\over{16\,\sqrt{15}} \\ {{32\,T^{10}-80\,T^8-80\,T^6+200\,T^4-150\,T^2-25
}}\over{80\,\sqrt{3}}\\{{\scriptstyle{T\,\left(32\,T^{10}-112\,T^8+48\,T^6+408\,T^4-
			150\,T^2+45\right)}}}\over{48\sqrt{15
}} \\{{\scriptstyle 64\,T^{12}-320\,T^{10}+464\,T^8+1120\,T^6-420\,T^4+1260\,T^2+
		315}}\over{96\,\sqrt{35}}\\\frac{T \left(64 T^{12}-448 T^{10}+1040 T^8+928 T^6-1284 T^4+3460 T^2+735\right)}{960}
\end{array}\right)\dfrac{-e^ {- T^2 }} {\sqrt{\pi}\left(2\,T^2+1\right)}.
$$



\bigskip


For $N=7$ we have the following expression for the corresponding commuting matrix $\mathscr{L}$


\begin{small}
	$$
	\left(
	\begin{array}{ccccccc}
	-\frac{2 T \left(14 T^2+27\right)}{\sqrt{3}} & 14 & 0 & 0 & 0 & 0 & 0 \\
14 & -\frac{5 T \left(22 T^2+21\right)}{4 \sqrt{3}} & \frac{5 \left(2 T^2+19\right)}{2 \sqrt{2}} & 5 \sqrt{\frac{5}{2}} T & \frac{5}{\sqrt{2}} & 0 & 0 \\
	0 & \frac{5 \left(2 T^2+19\right)}{2 \sqrt{2}} & -\frac{T \left(26 T^2+1\right)}{\sqrt{3}} & 2 \sqrt{\frac{5}{3}} \left(3 T^2+14\right) & 6 \sqrt{3} T & \sqrt{7} & 0 \\
	0 & 5 \sqrt{\frac{5}{2}} T & 2 \sqrt{\frac{5}{3}} \left(3 T^2+14\right) & -\frac{T \left(46 T^2-27\right)}{2 \sqrt{3}} & \frac{3}{2}  \sqrt{15} \left(2 T^2+5\right) & \frac{3 \sqrt{35} T}{2} & 1\\
	0 & \frac{5}{\sqrt{2}} & 6 \sqrt{3} T & \frac{3}{2} \sqrt{15} \left(2 T^2+5\right) & -\sqrt{3} T \left(6 T^2-5\right) & \sqrt{7} \left(5 T^2+7\right) & 2 \sqrt{5} T \\
	0 & 0 & \sqrt{7} & \frac{3 \sqrt{35} T}{2} & \sqrt{7} \left(5 T^2+7\right) & -\frac{T \left(42 T^2-29\right)}{4 \sqrt{3}} & \frac{\sqrt{35} \left(6 T^2+5\right)}{2 \sqrt{3}} \\
	0 & 0 & 0 & 1  & 2 \sqrt{5} T & \frac{\sqrt{35}\left(6 T^2+5\right)}{2 \sqrt{3}} & 0 \\
	\end{array}
	\right).
	$$

\end{small}

\bigskip

One may check that this matrix has simple spectrum.

\bigskip

Let us write $B_N$ for the truncated matrix  of size $N\times N$  made up to the first $N$ rows and columns of $B$ and $\Lambda_N$ the diagonal $N\times N$ matrix $\Lambda_N=\textrm{diag}(0,-6,-8,\ldots,-2(N+1))$. We write 
\begin{equation}\label{ads}
\textrm{ad}\ X(Y)=XY-YX 
\end{equation}
for the usual commutator and $\left(\textrm{ad}\ X\right)^n(Y)=\textrm{ad}\ X \left(\left(\textrm{ad}\ X\right)^{n-1}(Y)\right)$, $n\geq 2$.



\bigskip

\begin{remark}\label{combina}
The following relation holds true for every value of $N$, $N\geq 1$:
	
$$
-\dfrac{1}{40}\left(\textrm{ad}\  \Lambda_N \right)^4(B_N)+\left(\textrm{ad}\ \Lambda_N\right)^2(B_N)
-\dfrac{18}{5}B_N=0.
$$


	where $0$ is the zero matrix of size $N\times N$.
	
\end{remark}

A non obvious consequence of this is that if one tries to write the commuting matrix $\mathscr{L}$ in terms of $B_N$ and $\Lambda_N$ for  arbitrary values of $N$ and $T$,  in the spirit of Perline \cite{Per}, one would need to use monomials of degree higher than five. We will face a much better situation in the next section.

\begin{remark}

Interestingly, one can see that using the expression of $\Theta(x)$ in (\ref{coefs}) and the one of $L$ in (\ref{opLHermite}), one obtains
$$
-\dfrac{1}{40}\left(\textrm{ad}\  L \right)^4(\Theta)+\left(\textrm{ad}\ L\right)^2(\Theta)
-\dfrac{18}{5}\Theta=0.
$$
in the spirit of M. Reach in \cite{R1,R}, who was the first to consider mixed bispectral situations involving differential and difference operators.
\end{remark}


\bigskip

\begin{remark}
	We point out that we have made no use of the Darboux process which has always played an important role in terms of the bispectral property. For the case of exceptional Hermite polynomials discussed here we notice that the differential operator $L$ in (\ref{opLHermite})	
	 is Darboux connected, see \cite[Theorem 1.2 and Definition 3.7]{GGU2}, to the one for classical Hermite polynomials, namely$$
	\mathcal{H}=\dfrac{d^2}{dx^2}-2x\dfrac{d}{dx},
	$$
	by the relation	$$
L	\mathcal{S}=\mathcal{S}	\mathcal{H},
	$$
\end{remark}	
with
$$
\mathcal{S}=(2x^2+1)\dfrac{d^2}{dx^2}-4x
\dfrac{d}{dx}+4.
$$

\section{The exceptional Jacobi polynomials}\label{Jacobi}

\subsection{Bispectrality}

Let $\alpha$ and $\beta$ be the classical real Jacobi parameters, with $\alpha,\ \beta>-1$ and $\alpha \neq \beta$.

\bigskip

We write

\begin{equation}\label{parametrosJ}
 a=\frac{\beta-\alpha}{2}, b= \frac{\beta+\alpha}{(\beta-\alpha)}, c= b+\frac{1}{a}.
\end{equation}


We denote by $\hat{p}^{(\alpha,\beta)}_n(x)$ the sequence of exceptional Jacobi polynomials introduced in \cite{GU2}, where the set of missing degrees is equal to $\{0\}$. These polynomials are orthogonal in $(-1,1)$ with respect to the weight

 $$\rho(x) =\dfrac{(1-x)^{\alpha}(1+x)^{\beta}}{(x-b)^2}.$$ 

We use the expression
$$\hat{p}^{(\alpha,\beta)}_n(x)={{{\it b}\,{\it p^{(\alpha,\beta)}_{n-1}(x)-{\it p^{(\alpha,\beta)}_{n-2}(x)}}\over{2\,n+{\it \beta}+
 {\it \alpha}-2}}-{{ p^{(\alpha,\beta)}_{n-1}(x)\,\left(x-{\it b}\right)}\over{2}}}, \quad n\geq 1,$$


\bigskip

for the exceptional Jacobi polynomials $\hat{p}^{(\alpha,\beta)}_n$ in terms of the traditional Jacobi polynomials $p^{(\alpha,\beta)}_n(x)$.

It is very well known that the Jacobi polynomials are bispectral, i.e., they satisfy the recursion
(see for instance \cite[Chapter IV]{Sz})
$$c_n p^{(\alpha,\beta)}_{n-1}(x)+b_n p^{(\alpha,\beta)}_n(x) + a_n p^{(\alpha,\beta)}_{n+1}(x) = x p^{(\alpha,\beta)}_n(x), \quad n=0,1,2,\ldots,$$
with $a_n,b_n,c_n$ given below and $p^{(\alpha,\beta)}_{-1}=0,\ p^{(\alpha,\beta)}_{0}=1$,
as well as the differential equation
\begin{equation}
(1-x^2)(p^{(\alpha,\beta)}_n(x))^{\prime \prime} +(\beta-\alpha-(\alpha+\beta+2)x)(p^{(\alpha,\beta)}_n(x))^{\prime}(x)=-n(n+\alpha+\beta+1)p^{(\alpha,\beta)}_n(x), \quad n=0,1,2,\ldots.
\end{equation}

\bigskip

The values of $a_n,b_n,c_n$ above are given by $$a_n=\dfrac{2\,\left(n+1\right)\,\left(n+{\it \beta}+{\it \alpha}+1\right)}{
 \left(2\,n+{\it \beta}+{\it \alpha}+1\right)\,\left(2\,n+{\it \beta}+{\it \alpha}+
 2\right)}, \ b_n=\dfrac{\beta^2-{\it \alpha}^2}{\left(2\,n+{\it \beta}+{\it \alpha}\right)\,
 \left(2\,n+{\it \beta}+{\it \alpha}+2\right)}$$ and $c_n=\dfrac{2\,\left(n+{\it \alpha}\right)\,\left(n+{\it \beta}\right)}{\left(2
 \,n+{\it \beta}+{\it \alpha}\right)\,\left(2\,n+{\it \beta}+{\it \alpha}+1\right)
}.$

\bigskip




In this section we observe that the exceptional Jacobi polynomials $\hat{p}^{(\alpha,\beta)}_n(x)$ are also bispectral. They satisfy the five term recursion relation
\begin{equation}\label{relacionJac}
\hat{e}_n \hat{p}^{(\alpha,\beta)}_{n-2}(x)+\hat{d}_n \hat{p}^{(\alpha,\beta)}_{n-1}(x)+\hat{c}_n \hat{p}^{(\alpha,\beta)}_n(x) + \hat{b}_n \hat{p}^{(\alpha,\beta)}_{n+1}(x)+ \hat{a}_n \hat{p}^{(\alpha,\beta)}_{n+2}(x) = (x-b)^2 \hat{p}^{(\alpha,\beta)}_n(x),
\end{equation}
 where the explicit expressions of the entries $\hat{e}_n$, $\hat{d}_n$, $\hat{c}_n$, $\hat{b}_n$ and $\hat{a}_n$ are given below, as well as the differential equation (see \cite[Section 2]{GU2})
\begin{equation}\label{Jacdifeq}
L_{\alpha ,\beta}\hat{p}^{(\alpha,\beta)}_n(x)=(n-1)(\alpha+\beta+n)\hat{p}^{(\alpha,\beta)}_n(x), \quad n=1,2\ldots,
\end{equation}
where the operator $L_{\alpha ,\beta}$ is given by
\begin{equation}\label{difopJacobi}
L_{\alpha ,\beta}=(x^2-1)\frac{d^2}{dx^2}+2a\left(\frac{1-bx}{b-x}\right)\left((x-c) \frac{d}{dx}-1 \right),
\end{equation}
 with the parameters $a$, $b$ and $c$ given in (\ref{parametrosJ}).
\bigskip

\bigskip

The existence of a five term recursion relation was established in \cite{STZ}. This can be expressed by saying that the vector
\begin{equation}\label{Psi}
\Psi(x)=(\hat{p}^{(\alpha,\beta)}_{1}(x),\hat{p}^{(\alpha,\beta)}_{2}(x),\hat{p}^{(\alpha,\beta)}_{3}(x),\ldots)^t
\end{equation}

satisfies the relation $$\hat{B} \Psi(x)  = (x-b)^2 \Psi(x),$$

where $\hat{B}$ is the pentadiagonal  matrix where  its $n-th$ row is given by $$(...0,0,0,\hat{e}_n,\hat{d}_n,\hat{c}_n,\hat{b}_n,\hat{a}_n,0,0,...).$$

\bigskip

The matrix $\hat{B}$ is not symmetric, but it can be symmetrized into a matrix $B$ with
$B= \mathcal{D} \hat{B} \mathcal{D}^{-1}$, where the diagonal matrix $\mathcal{D}$ with entries $\mu_i$, $i=1,2,\ldots$, gives rise to a sequence of
orthonormal exceptional Jacobi polynomials, that written in terms of a  vector $\widetilde{\Psi}(x)$ as in (\ref{Psi}), reads as follows: 

\begin{equation}\label{thetax}
B  \widetilde{\Psi}(x) = (x-b)^2 \widetilde{\Psi}(x).
\end{equation}

\bigskip

The entries of the recurrence relation (\ref{relacionJac}) are given as follows:

\bigskip


$$\hat{a}_n={{4\,n\,\left(n+1\right)\,\left(n+{\it \beta}+{\it \alpha}\right)\,\left(n
 +{\it \beta}+{\it \alpha}+1\right)}\over{\left(2\,n+{\it \beta}+{\it \alpha}-1
 \right)\,\left(2\,n+{\it \beta}+{\it \alpha}\right)\,\left(2\,n+{\it \beta}+
 {\it \alpha}+1\right)\,\left(2\,n+{\it \beta}+{\it \alpha}+2\right)}},$$

 \bigskip


$$\hat{b}_n=-{{16n\,\left({\it \beta}+{\it \alpha}\right)\left(n+{\it \alpha}\right)\,
 \left(n+{\it \beta}\right)\,\left(n+{\it \beta}+{\it \alpha}\right)}\over{
 \left({\it \beta}-{\it \alpha}\right)\,\left(2\,n+{\it \beta}+{\it \alpha}-2
 \right)\,\left(2\,n+{\it \beta}+{\it \alpha}-1\right)\,\left(\alpha+{\it \beta}+2n\right)\left(\alpha+{\it \beta}+2n+2\right)}},
$$

\bigskip


\begin{eqnarray*}
\hat{c}_n&=&-{{2\,n\,\left(n+1\right)\,\left(n+{\it \beta}-1\right)\,\left(n+
 {\it \beta}+2\right)}\over{2\,n+{\it \beta}+{\it \alpha}+1}}+{{4\,n^2\,\left(n+{\it \beta}-1\right)\,\left(n+{\it \beta}+1\right)
 }\over{2\,n+{\it \beta}+{\it \alpha}}}-{{4\,\left(n-1\right)^2\,\left(n+{\it \beta}-2\right)\,\left(n+
 {\it \beta}\right)}\over{2\,n+{\it \beta}+{\it \alpha}-2}}\\
&+&{{2\,\left(n-2\right)\,\left(n-1\right)\,\left(n+{\it \beta}-3\right)
 \,\left(n+{\it \beta}\right)}\over{2\,n+{\it \beta}+{\it \alpha}-3}}+{{4\,{\it \beta}^2}\over{\left({\it \beta}-{\it \alpha}\right)^2}}-{{4\,
 {\it \beta}}\over{{\it \beta}-{\it \alpha}}},
\end{eqnarray*}

$$\hat{d}_n=-{{16\,\left({\it \beta}+{\it \alpha}\right)\,\left(n+{\it \alpha}-2\right)\,
		\left(n+{\it \alpha}\right)\,\left(n+{\it \beta}-2\right)\,\left(n+{\it \beta}
		\right)}\over{\left({\it \beta}-{\it \alpha}\right)\,\left(2\,n+{\it \beta}+
		{\it \alpha}-4\right)\,\left(2\,n+{\it \beta}+{\it \alpha}-2\right)\,\left(2\,n
		+{\it \beta}+{\it \alpha}-1\right)\,\left(2\,n+{\it \beta}+{\it \alpha}\right)}},$$

$$\hat{e}_n={{4\,\left(n+{\it \alpha}-3\right)\,\left(n+{\it \alpha}\right)\,\left(n+
		{\it \beta}-3\right)\,\left(n+{\it \beta}\right)}\over{\left(2\,n+{\it \beta}
		+{\it \alpha}-4\right)\,\left(2\,n+{\it \beta}+{\it \alpha}-3\right)\,\left(2\,
		n+{\it \beta}+{\it \alpha}-2\right)\,\left(2\,n+{\it \beta}+{\it \alpha}-1\right)
}}.$$


\bigskip


Finally,
the expression for the ratio $\left(\dfrac{\mu_{n+1}}{\mu_n}\right)^2$ is given by
$${{n\,\left(n+{\it \alpha}\right)\,\left(n+{\it \beta}\right)\,\left(n+
 {\it \beta}+{\it \alpha}\right)\,\left(2\,n+{\it \beta}+{\it \alpha}+1\right)
 }\over{\left(n+{\it \alpha}-1\right)\,\left(n+{\it \alpha}+1\right)\,\left(n
 +{\it \beta}-1\right)\,\left(n+{\it \beta}+1\right)\,\left(2\,n+{\it \beta}+
 {\it \alpha}-1\right)}}$$

 \bigskip

 and this determines the diagonal matrix $\mathcal{D}$ with positive coefficients up to a scalar.

\bigskip

We have by now identified the differential operator $L$ and the difference operator $B$ as well as 
the diagonal operators $\Lambda(k)=\Lambda(n)$ and $\Theta(x)=(x-b)^2$ (see (\ref{bispL}) and (\ref{bispB})) that give us a bispectral situation for a sequence of orthonormal exceptional Jacobi polynomials.


\subsection{The commuting matrix}


We find a narrow banded matrix $\mathscr{L}$ that commutes with the band-time-band limiting matrix of inner products. This is illustrated in the case of size $N=7$, the band limiting parameter $T=1/3$ and the Jacobi parameters $\alpha=3$ and $\beta=4$:

\bigskip


$$
\left(
\begin{array}{ccccccc}
	-\frac{316897}{2544 \sqrt{4290}} & \frac{26486 \sqrt{\frac{2}{429}}}{583} & -\frac{969 \sqrt{\frac{7}{1430}}}{1166} & 0 & 0 & 0 & 0 \\
	\frac{26486 \sqrt{\frac{2}{429}}}{583} & -\frac{200311 \sqrt{\frac{5}{858}}}{6996} & \frac{82042 \sqrt{\frac{14}{143}}}{7579} & -\frac{5491}{7579 \sqrt{143}} & 0 & 0 & 0 \\
	-\frac{969 \sqrt{\frac{7}{1430}}}{1166} & \frac{82042 \sqrt{\frac{14}{143}}}{7579} & -\frac{21492811}{159159 \sqrt{4290}} & \frac{2081089 \sqrt{\frac{7}{2145}}}{37895} & -\frac{148257}{385840 \sqrt{65}} & 0 & 0 \\
	0 & -\frac{5491}{7579 \sqrt{143}} & \frac{2081089 \sqrt{\frac{7}{2145}}}{37895} & -\frac{24836371}{220480 \sqrt{4290}} & \frac{19323 \sqrt{\frac{14}{65}}}{3445} & -\frac{3249 \sqrt{\frac{3}{130}}}{16960} & 0 \\
	0 & 0 & -\frac{148257}{385840 \sqrt{65}} & \frac{19323 \sqrt{\frac{14}{65}}}{3445} & -\frac{7219271}{89040 \sqrt{4290}} & \frac{323 \sqrt{\frac{154}{65}}}{265} & -\frac{19 \sqrt{\frac{7}{3}}}{2544} \\
	0 & 0 & 0 & -\frac{3249 \sqrt{\frac{3}{130}}}{16960} & \frac{323 \sqrt{\frac{154}{65}}}{265} & -\frac{437291}{10176 \sqrt{4290}} & 1 \\
	0 & 0 & 0 & 0 & -\frac{19 \sqrt{\frac{7}{3}}}{2544} & 1 & 0 \\
\end{array}
\right).
$$

  The commuting matrix $\mathscr{L}$ of size $N$ is unique once we ask for $\mathscr{L}_{N,N}=0$ and $\mathscr{L}_{N,N-1}=1$.

\bigskip

As in the previous section, one writes  $B_N$ for the truncated matrix  of size $N\times N$  made up to the first $N$ rows and columns of $B$ and $\Lambda_N$ the diagonal $N\times N$  matrix  whose entry in the position $(n,n)$ is equal to the eigenvalue   $\lambda_n=(n-1)(n+\alpha+\beta)$ in the differential equation (\ref{Jacdifeq}).

\bigskip

Attempting to express the matrix $\mathscr{L}$ above as a linear combination of very simple monomials in terms of $B_N$ and $\Lambda_N$, with $N=7$, one obtains 


\begin{eqnarray} \label{PerlineL}
\mathscr{L}&=&\gamma_1I+\gamma_2\Lambda_N+\gamma_3B_N+\gamma_4\Lambda_N^2+\gamma_5\left(\Lambda_N B_N+ B_N \Lambda_N\right)+\gamma_6\left(  B_N \Lambda_N^2+\Lambda_N ^2 B_N\right)\\&+& \gamma_7\Lambda_N B_N \Lambda_N +
\gamma_8\left[\left(\Lambda_N^3 B_N + B_N \Lambda_N^3\right)-\left(\Lambda_N B_N \Lambda_N^2+\Lambda_N^2 B_N \Lambda_N\right)\right], \nonumber
\end{eqnarray}
where
$${\it \gamma_1}={{7387129}\over{53\,2^{{{5}\over{2}}}\,3^{{{3}\over{2}}}
		\,5^{{{3}\over{2}}}\,\sqrt{11}\,\sqrt{13}}},\ {\it \gamma_2}=-{{4796873}\over{371\,2^{{{9}\over{2}}}\,3^{{{5}\over{2}}
 }\,5^{{{3}\over{2}}}\,\sqrt{11}\,\sqrt{13}}}, \ {\it \gamma_3}=-{{116603\,\sqrt{3}}\over{53\,2^{{{15}\over{2}}}\,\sqrt{5
 }\,\sqrt{11}\,\sqrt{13}}},$$
$$ {\it \gamma_4}=-{{323\,\sqrt{11}}\over{371\,2^{{{9}\over{2}}}\,3^{{{5
 }\over{2}}}\,\sqrt{5}\,\sqrt{13}}},
\ {\it \gamma_5}={{7429\,3^{{{5}\over{2}}}}\over{371\,2^{{{19}\over{2}}}\,
 5^{{{3}\over{2}}}\,\sqrt{11}\,\sqrt{13}}},\ {\it \gamma_6}={{77843}\over{371\,2^{{{15}\over{2}}}\,3^{{{5}\over{2}}}
 \,5^{{{3}\over{2}}}\,\sqrt{11}\,\sqrt{13}}}, \ {\it \gamma_7}=-{{305881}\over{371\,2^{{{17}\over{2}}}\,3^{{{5}\over{2}}
 }\,5^{{{3}\over{2}}}\,\sqrt{11}\,\sqrt{13}}},$$ and
$${\it \gamma_8}=-{{323}\over{371\,2^{{{19}\over{2}}}\,3^{{{5}\over{2}}}\,
 \sqrt{5}\,\sqrt{11}\,\sqrt{13}}}.$$

\bigskip

This expression, of course, can be written in terms of commutators and anti-commutators as in Remark \ref{combina}.

\bigskip

It is worth pointing out that a  similar combination, with the same monomials as above, holds true for arbitrary size $N\geq7$. This yields a non trivial extension of the results in \cite{Per}.






\subsection{The commuting differential operator}

For given values of $N$ and $T$, consider the integral kernel
$$
K_N(x,y)=\sum_{i=1}^N\dfrac{ \hat{p}^{(\alpha,\beta)}_i(x) \hat{p}^{(\alpha,\beta)}_i(y)}{ ||\hat{p}^{(\alpha,\beta)}_i  ||^2 }
$$
acting on $L^2(-1,T)$.

\bigskip

In perfect agreement with the expression of $\mathscr{L}$ given in (\ref{PerlineL}) one can build a differential operator commuting with the integral operator with the kernel given above, by taking an appropriate linear combination in terms of the operators $L=L_{\alpha,\beta}$ and $\Theta(x)=(x-b)^2$ (see (\ref{difopJacobi}) and (\ref{thetax})). The linear combination involves the following operators:
$$ 
L,\ \Theta,\ L^2,\ L\Theta+\Theta L,\ L^2\Theta+\Theta L^2, \ L\Theta L,\ (L\Theta L^2+L^2\Theta L)-(L^3\Theta+\Theta L^3).
$$

This operator can be written in a more explicit form, namely
\begin{equation}\label{superOpJacobi}
\dfrac{1}{\rho(x)}\left( \dfrac{d^2}{dx^2}\mathscr{A}(x)\dfrac{d^2}{dx^2}+\dfrac{d}{dx}\mathscr{B}(x)\dfrac{d}{dx}+\mathscr{C}(x)\right)
\end{equation}
where $$\rho(x) =\dfrac{(1-x)^{\alpha}(1+x)^{\beta}}{(x-b)^2}$$ is the orthogonality weight of the exceptional Jacobi polynomials $\hat{p}^{(\alpha,\beta)}_n(x)$  displayed above and the functions $\mathscr{A}(x)$, $\mathscr{B}(x)$, $\mathscr{C}(x)$ have the form
\begin{eqnarray*}
\mathscr{A}(x)&=&\dfrac{(x-1)^{\alpha+2}(x+1)^{\beta+2}(x-T)^2}{(x-b)^2},
\quad 
\mathscr{B}(x)=\dfrac{(x-1)^{\alpha+1}(x+1)^{\beta+1}(x-T)}{(x-b)^4}\mathcal{P}_3(x),\\
\mathscr{C}(x)&=&\dfrac{(x-1)^{\alpha}(x+1)^{\beta}}{(x-b)^5}\mathcal{P}_5(x).
\end{eqnarray*}

\bigskip

Notice that the differential operator given above in (\ref{superOpJacobi}) is selfadjoint in $L^2((-1,T),\rho(x)dx)$.

\bigskip

Here $\mathcal{P}_3(x)$ and $\mathcal{P}_5(x)$ are polynomials in $x$ of degree three and five respectively.

\medskip

One can describe the dependence of $\mathcal{P}_3(x)$ and $\mathcal{P}_5(x)$ on the parameters $T$ and $N$ in more detail, namely
$$
\mathcal{P}_3(x)=q_3(N)x^3+(q_2(N)+\sigma_2 T)x^2+(q_1(N)+\sigma_1 T)x+q_0(N)+\sigma_0 T,
$$
with $q_i$, $i=0,\ldots,3$, quadratic polynomials in $N$ and $\sigma_i$, $i=0,1,2$, constants. Moreover, $q_3$ is up to a multiplicative constant 
$(N-1)(N+\alpha+\beta+1)$.

\medskip

In the case of $\mathcal{P}_5(x)$ one gets
\begin{eqnarray*}
	\mathcal{P}_5(x)&=&\delta_5(N)x^5+(\delta_4(N)+\lambda_4(N)T)x^4+(\delta_2(N)+\lambda_2(N)T+\gamma_2T^2)x^2+(\delta_1(N)+\lambda_1(N)T+\gamma_1T^2)x\\
	&+&\delta_0(N)+\lambda_0(N)T+\gamma_0T^2.
\end{eqnarray*}

By adding to $\mathscr{C}(x)$ a constant 
one can make, as above, one of the coefficients of $\mathcal{P}_5(x)$ vanish.  Here $\delta_5(N)$ is up to a multiplicative constant equal to $(N-1)N(N+\alpha+\beta)(N+\alpha+\beta+1)$, the coefficients $\delta_i$, $i=0,\ldots 4$, are polynomials in  $N$ of degree four (not as nice as $\delta_5$), $\lambda_i$, $i=0,1,2,4$, are polynomials in $N$ of degree two, and $\gamma_i$, $i=0,1,2$, are constants independent of $N$.

%


\bigskip

\begin{remark}
 Using the notation in (\ref{ads}) one can see that the ad conditions in \cite{R1,R} expressing the bispectral property take the following form in this case
 \begin{eqnarray*}
 &&\left(\textrm{ad}\ L_1\right)^5(\Theta)-10\ \left(\textrm{ad}\ L_1\right)^4(\Theta)+\left(\textrm{ad}\ L_1\right)^3\left(\Theta\right)\left(33-20L_1\right)+\left(\textrm{ad}\ L_1\right)^2\left(\Theta\right)\left(-40+64L_1\right)\\&&+16\ \textrm{ad}\ L_1\left(\Theta\right)\left(\mathcal{I}-5L_1-4L_1^2\right) =0,
  \end{eqnarray*}
 	where $\mathcal{I}$ is the identity operator and $L_1=L-\dfrac{(\alpha+\beta+1)^2}{4}\mathcal{I}$ is a proper shift of $L$. The operator $\Theta(x)=(x-b)^2$ is the same as above.
 	
 	\medskip
 	
 	This more complicated form of the ad conditions in \cite{R1,R} arises because the eigenvalue in the Jacobi case is a quadratic function of $n$.

\end{remark}

\bigskip

\begin{remark}
For the case of exceptional Jacobi polynomials discussed here we notice that the operator 	$L_{\alpha,\beta}$ in (\ref{difopJacobi})	 is Darboux connected to the one for classical Jacobi polynomials, namely$$
J_{\alpha,\beta}=(x^2-1)\dfrac{d^2}{dx^2}-(\beta-\alpha-(\alpha+\beta+2)x)\dfrac{d}{dx},
$$
by the relation	$$
	L_{\alpha,\beta}\mathcal{S}=\mathcal{S}J_{\alpha,\beta},
	$$
\end{remark}	
with
$$
\mathcal{S}=\mathcal{S}_2(x)\dfrac{d^2}{dx^2}+\mathcal{S}_1(x)
\dfrac{d}{dx}+\mathcal{S}_0(x),
$$
where \begin{eqnarray*}
\mathcal{S}_2(x)&=&(x^2-1)((\beta-\alpha)x-\beta-\alpha), \\
\mathcal{S}_1(x)&=&(\beta-\alpha)(\alpha+\beta+1)x^2-2(\beta(\beta+1)+\alpha(\alpha+1))x+(\beta-\alpha)(\alpha+\beta+1),\\
\mathcal{S}_0(x)&=& \alpha\beta ((\beta - \alpha) x - \beta - \alpha - 2).
\end{eqnarray*}

\section{The exceptional Laguerre polynomials}\label{Laguerre}


Let $\alpha>0$ and $L^{(\alpha)}_n(x)$, $n\geq 0$, denote the classical Laguerre polynomials orthogonal with respect
to the weight $e^{-x}x^\alpha$ in the interval $(0,+\infty)$. We consider the sequence of exceptional Laguerre polynomials $\hat{L}^{(\alpha)}_n(x)$ introduced in \cite{GU2}, orthogonal with respect to the weight

$$\frac{e^{-x}x^\alpha}{(x+\alpha)^2} $$

where the set of missing degrees is $\{0\}$.

\medskip

These polynomials can be expressed in terms of the classical Laguerre polynomials by the relation (see \cite[section 6.2]{GU2})
$$\hat{L}^{(\alpha)}_n(x)=-(x+\alpha+1)L_{n-1}^{(\alpha)}+L_{n-2}^{(\alpha)},\quad n=1,2,\ldots .$$

\bigskip

The squared norms of these polynomials are given by
$$||\hat{L}^{(\alpha)}_n(x)||^2=\frac{(\alpha+n)\Gamma(n+\alpha)}{(\alpha+n-1)(n-1)!}. $$
This expression is essentially  in \cite[section 6.2]{GU2}). We write $\left( \widetilde{L}^{(\alpha)}_n(x)\right)_{n\geq 1}$ for the sequence of orthonormal polynomials.

\bigskip


Here the relevant operators $L=L_{\alpha}$ and $\Theta(x)=\Theta_{\alpha}(x)$ (see (\ref{bispL}) and (\ref{bispB})) are given by $\Theta_{\alpha}=(x+\alpha)^2$ and
\begin{equation}\label{Laguerre_op}
L_{\alpha}=-x\dfrac{d^2}{dx^2}+\left(\dfrac{x-\alpha}{x+\alpha}\right)\left((x+\alpha+1)\dfrac{d}{dx}-\mathcal{I}\right).
\end{equation}

\bigskip

\begin{remark}
In terms of $L$ and $\Theta$ the ad conditions in \cite{R1,R} take, for this case, the simple form (see (\ref{ads}) for the notation)
$$
\left(\textrm{ad}\ L\right)^5(\Theta)-5\left(\textrm{ad}\ L\right)^3(\Theta)+4\ \textrm{ad}\ L(\Theta)=0.
$$

\end{remark}

\bigskip

Turning our attention to the commuting property, we define the $N\times N$ matrix $M$ of truncated inner products depending on a real parameter $T$, whose entries are given by:

$$M_{m,n}=\int_{0}^{T}\widetilde{L}^{(\alpha)}_m(x)\widetilde{L}^{(\alpha)}_n(x)\frac{e^{-x}x^\alpha}{(x+\alpha)^2} dx,\quad n,m=1,2,\ldots, N. $$
For each fixed value of the parameters $N$ and $T$ one looks for a ``narrow banded" commuting matrix $\mathscr{L}$.

\bigskip

We exhibit the commuting matrix $\mathscr{L}$ of size $N=7$, for the special choice of the Laguerre  parameter $\alpha=7$.


\bigskip

\begin{scriptsize}
$$
\begin{pmatrix}
 -{{\sqrt{7}\,\left(3\,T^2+16\,T-260\right)}}\over{110\sqrt{2}} & -{{6\,\left(T+19\right)}\over{55}} & {{\sqrt{5}

	}\over{11\,\sqrt{2}}} & 0 & 0 & 0 &0 \\
-{{6\,\left(T+19\right)}\over{55}} & -{{11\,T^2+32\,T-1232

	}\over{66\,\sqrt{14}}} & -{{2\,\sqrt{5}\,\left(T+

		17\right)}\over{11\,\sqrt{7}}} &{2\sqrt{2}}\over{\sqrt{3}
	\,\sqrt{7}\,\sqrt{11}} & 0 &0 & 0\\
	 {{\sqrt{5}}\over{11\,\sqrt{2}}} & -{{2\,\sqrt{5}\,\left(T+17

			\right)}\over{11\,\sqrt{7}}} & -{{23\,T^2+8\,T-3060}\over{165\,

			\sqrt{14}}} & -{{2\,\sqrt{3}\,\left(T+15\right)}\over{

			\sqrt{5}\,\sqrt{7}\,\sqrt{11}}} & {18\sqrt{2}}\over{55\,

			\sqrt{7}} & 0 &0 \\
		 0 & {2\sqrt{2}}\over{\sqrt{3}\,\sqrt{7}\,\sqrt{11}
	} & -{{2\,\sqrt{3}\,\left(T+15\right)}\over{\sqrt{5}\,\sqrt{7}\,

				\sqrt{11}}} & -{{6\,T^2-14\,T-893}\over{55\,\sqrt{14}}} & -

		{{2\,\sqrt{3}\,\left(T+13\right)}\over{\sqrt{5}\,\sqrt{7}\,\sqrt{11}

		}} & {{\sqrt{13}}\over{11\,\sqrt{7}}} &0\\
0 &0 &{18\sqrt{2}}\over{55\,\sqrt{7}} & -{{2\,

		\sqrt{3}\,\left(T+13\right)}\over{\sqrt{5}\,\sqrt{7}\,\sqrt{11}}} & 

-{{25\,T^2-128\,T-3964}\over{330\sqrt{14}}} & -

{{\sqrt{26}\left(T+11\right)}\over{\sqrt{3}\,\sqrt{5}\,

		\sqrt{7}\,\sqrt{11}}} &{1}\over{\sqrt{330}} \\
		0 &0 & 0 &{{\sqrt{13}}\over{11\,\sqrt{7}}} & -{{\sqrt{26}

			\left(T+11\right)}\over{\sqrt{3}\,\sqrt{5}\,\sqrt{7}\,

			\sqrt{11}}} & -{{13\,T^2-104\,T-2100}\over{330\sqrt{14}
			}} & -{{\sqrt{26}\,\left(T+9\right)}\over{55}}
\\
		0 & 0 & 0 &0 & {1}\over{\sqrt{330}}& -{{\sqrt{26}\left(T+9\right)}\over{55}} &

		0 
		\end{pmatrix}.
 $$
\end{scriptsize}

\bigskip

As before, with a proper normalization this matrix is unique.

\bigskip

\begin{remark}
For the case of exceptional Laguerre polynomials discussed here we notice that the differential operator 	$L_{\alpha}$ in (\ref{Laguerre_op})
	 is Darboux connected to the one for classical Laguerre polynomials, namely$$
	\mathcal{L}_{\alpha}=-x\dfrac{d^2}{dx^2}-(\alpha+1-x)\dfrac{d}{dx},
	$$
	by the relation	$$
		L_{\alpha}		\mathcal{S}=\mathcal{S}\mathcal{L}_{\alpha},
	$$
	with
	$$ S=x(x+\alpha)\dfrac{d^2}{dx^2}-(x^2-\alpha^2-\alpha)\dfrac{d}{dx}-\alpha(x+\alpha+1).$$
\end{remark}

\section{The benefit of having a commuting local matrix}\label{numeric}


The previous sections have exhibited, for a collection of examples of exceptional orthogonal polynomials, a pair of matrices. The first one, a full matrix $M=M_{T,N}$, is obtained by forming the inner products of the normalized OP over a restricted range in physical space, thus implementing ``time limiting" with parameter $T$. In the resulting $N \times N$ matrix the parameter $N$ implements ``band limiting". In each case the second matrix, denoted by $\mathscr{L}=\mathscr{L}_{T,N}$, is a narrow banded one that commutes with the first matrix, and has simple spectrum. From a numerical point of view the entire purpose of the search for this second matrix is that it reduces the problem of computing the eigenvectors of the first one, a seriously ill-conditioned one, into a very well conditioned one.

	\bigskip

The eigenvectors of $M_{T,N}$ are of paramount importance since they give the singular vectors of the signal processing problem at hand, as described in Section \ref{intro}.

	\bigskip

We display below the results of some small size numerical computations
that illustrate the problem of computing the eigenvectors of a full matrix, such as $M_{T,N}$, some of whose eigenvalues are very close together. In each case, we give the eigenvalues of both the full matrix  $M_{T,N}$ and those of the narrow banded matrix $\mathscr{L}_{T,N}$.  It should be clear that in the case when $N$ is large the problems indicated below get to be much worse. Our point is that they  already appear for small values of $N$.

	\bigskip

We use the QR algorithm as implemented in LAPACK. In each of the three situations, Hermite, Jacobi and Laguerre, we will denote	by $X_M$ the matrix of eigenvectors  of $M_{T,N}$ (normalized and given as columns of $X_M$. We will denote by $Y_L$ the matrix of eigenvectors  of $\mathscr{L}_{T,N}$ (normalized and given as columns of $Y_L$).

	\bigskip

In theory the eigenvectors of $M_{T,N}$  should agree (up to order and signs) with those of $\mathscr{L}_{T,N}$.
 If we compute the matrix
	of inner products given by $$Y_L^T X_M$$ we expect to have the identity matrix up to some permutation and possibly some signs due to the normalization of the eigenvectors which are the columns of $X_M$ and $Y_L$.


	\subsection{Hermite}

	The choice of parameters  is  $N=7$ and $T=5$.
	
	\bigskip

	The eigenvalues of $M_{T,N}$
are
$$
1.,1.,1.,1.,1.,1.,0.999989,
$$

	and those of $\mathscr{L}_{T,N}$
are
	$$
	-2186.9,-2127.47,-1985.6,-1690.49,-1227.68,-611.685,100.033.
	$$

	
	\bigskip
	
	The matrix $Y_L^T X_M$ is
	
			\begin{scriptsize}
$$
\left(
\begin{array}{ccccccc}
	0.0663098 & 0.384538 & 0.914298 & -0.108586 & 0.00078648 & 0.0000796039 & -0.0000275025 \\
	0.250947 & 0.784559 & -0.396329 & -0.405483 & -0.000552959 & -0.0000929769 & -0.0000926722 \\
	-0.568584 & 0.481897 & -0.0828753 & 0.661528 & 0.000385822 & 0.000140271 & 0.000264937 \\
	0.780602 & 0.0661256 & -0.0106201 & 0.621429 & -0.000970967 & -0.000161057 & -0.000249407 \\
	0.00106371 & 9.69168\times 10^{-6} & -0.000916543 & 0.000209371 & 0.999999 & -0.00028544 & -0.000399834 \\
	0.000223739 & -0.000014601 & -0.0000999801 & -0.0000216885 & 0.000285012 & 1. & -0.000257904 \\
	0.000370889 & -0.0000278978 & 7.33262\times 10^{-6} & -0.0000607604 & 0.000399533 & 0.000257706 & 1. \\
\end{array}
\right)
$$

	\end{scriptsize}

	\subsection {Jacobi}

	The choice of parameters is given by $\alpha=4$, $\beta=3$,  $N=7$ and $T=1/3$.

	\bigskip

	The eigenvalues of $M_{T,N}$ are
$$ 0.00271069 , 0.0397751 , 0.752568 , 0.977212 , 1.0016 , 1.0
	, 1.0.$$ 
	Notice that the value 1.0016 above is due to numerical instability. Those of $\mathscr{L}_{T,N}$ are
	$$-225.02 , -131.013 , -47.8296 , 180.78 , 22.1365 , 119.864	, 74.8259.$$

	\bigskip

The matrix $Y_L^T X_M$ is

\begin{scriptsize}
$$	
\begin{pmatrix}	
 0.00747678 & -0.00879137 & -0.00787297 & 0.997507 & 

	0.00605771 & 0.0686547 & -0.00573508 \\
 0.0299278 & -0.0351898 & -0.0315137 & 0.239956 & 0.0242519

	& 0.968588 & -0.0232993 \\
 0.188064 & -0.221129 & -0.198037 & 0.0204167 & 0.154553 & -

	0.0312631 & 0.922623 \\ 0.470392 & -0.553098 & -0.49561 & 0.00362164 & 0.476008 & -0.00544977 & 0.0235557 \\ 
	-0.533333 & 0.6271 & 0.556811 & 2.82016 \times 10^{-4} & -
		0.110716 & -4.23774 \times 10^{-4} & 0.00171454 \\
		0.643064 & -0.753535 & -0.136571 & -7.37868 \times 10^{-7}
			& 3.89179 \times 10^{-4} & 1.10886 \times 10^{-6} & -
			4.50479 \times 10^{-6}\\
			-0.994134 & 0.108155 & 9.28542 \times 10^{-5} & 
			4.00526 \times 10^{-10} & -2.1141 \times 10^{-7} & -
			6.01908 \times 10^{-10} & 2.44529 \times 10^{-9} \\
\end{pmatrix}
	$$
\end{scriptsize}

	\subsection{Laguerre}

	The choice of parameters is given by $\alpha=7$, $N=7$ and $T=1/2$.

	\bigskip

	The eigenvalues of $M_{T,N}$ are


	$$
	1.0856\times 10^{-4} , 6.4646 \times 10^{-9} , 
	2.4992 \times 10^{-13} , -7.9725 \times 10^{-17} , -
	1.8132 \times 10^{-18} , 9.8169 \times 10^{-17} ,
		 6.469 \times 10^{-17} ,
$$
	 and those of $\mathscr{L}_{T,N}$ are
	 
	  $$169.43 , 119.107 , 78.6334 , 46.531 , 21.5219 , 2.5574 , -11.0609 $$


	\bigskip
	
	The matrix $Y_L^T X_M$ is

\begin{scriptsize}
	$$
	\begin{pmatrix}	
 4.53595 \times 10^{-14} & 7.10889 \times 10^{-15} & -
	3.33066 \times 10^{-15} & 4.84889 \times 10^{-14} & 
	4.16888 \times 10^{-14} & -4.76174 \times 10^{-13} & -1.0  \\
 -2.24211 \times 10^{-9} & 2.66194 \times 10^{-10} & 
1.39646 \times 10^{-9} & -1.14752 \times 10^{-9} & -
	1.08174 \times 10^{-8} & 0.999999 & -4.76896 \times 10^{-13}	\\
 6.62053 \times 10^{-5} & -1.75441 \times 10^{-5} & -
6.23886 \times 10^{-5} & -2.52534 \times 10^{-4} & 0.999999 & 
	1.08165 \times 10^{-8} & 4.19664 \times 10^{-14} \\
0.218337 & 0.447717 & 0.653743 & 0.569647 & 
1.78046 \times 10^{-4} & 1.13123 \times 10^{-10} & 
	3.85524 \times 10^{-14} \\
0.423228 & 0.637769 & -0.0176169 & -0.643286 & -
1.80397 \times 10^{-4} & 6.32598 \times 10^{-11} & -
7.10542 \times 10^{-15}\\
-0.104047 & 0.513037 & -0.717294 & 0.459842 & 
8.72632 \times 10^{-5} & 1.1602 \times 10^{-9} & 
2.35263 \times 10^{-14}	\\
	 -0.873152 & 0.359929 & 0.240438 & -0.224153 & 
	2.24938 \times 10^{-5} & -2.64713 \times 10^{-9} & -
	4.8397 \times 10^{-14} 
	\end{pmatrix}
$$
\end{scriptsize}

	
	\subsection{Conclusions}

	Observe that some of the entries of these matrices $Y^T_LX_M$ are indeed very close to the theoretically correct values, while	others are terribly off. The reason	is that a few eigenvalues of the full matrix of inner products $M_{T,N}$ are just too close together. This produces numerical instability in the computation of the corresponding eigenvectors. On the other hand, all the eigenvalues of the commuting matrix $\mathscr{L}_{T,N}$ are nicely separated and the corresponding eigenvectors can be trusted.


In summary, a good way to obtain reliable numerical values for the
eigenvectors of the global matrix $M$ is to forget about $M$ altogether	and to compute numerically the eigenvectors of $\mathscr{L}$. Not only we will then be dealing with a very sparse matrix for which the QR algorithm works	very fast (most of the work is avoided) but the problem	is numerically very well conditioned.

The main point of this paper is to show that the miracle of Slepian, Landau and Pollak, which we described in section \ref{intro}, applies to the examples of exceptional orthogonal polynomials discussed here and has important consequences.



\bigskip

\section*{Funding}

The work of the first author was  partially supported by  PID2021-124332NB-C21 (FEDER(EU)/Ministerio de Ciencia e Innovaci\'on-Agencia Estatal de Investigaci\'on) and FQM-262 (Junta	de Andalucía).

\section*{Declarations}


{\bf Conflict of interest} The authors declared that they have no conflict of interest.

\end{document}